\documentclass[12pt,leqno,twoside]{amsart}
\usepackage{amssymb,amsfonts,amsmath,amsthm}
\usepackage{latexsym}
\usepackage{verbatim}
\usepackage{epsfig}

\usepackage[T1]{fontenc}
\usepackage[latin1]{inputenc}

\newtheorem{theorem}{Theorem}[section]

\newtheorem{proposition}[theorem]{Proposition}
\theoremstyle{remark}
\newtheorem{remark}[theorem]{\sc Remark}
\theoremstyle{remark}

\theoremstyle{definition}
\newtheorem{definition}[theorem]{Definition}
\theoremstyle{remark}

\theoremstyle{remark}

\theoremstyle{remark}

\numberwithin{equation}{section}

\renewcommand{\Box}{_\square}    

\newcommand{\cal}{\mathcal}

\newcommand{\reg}{{\rm{reg}}}

\newcommand{\Sing}{\mathop{{\rm{Sing}}}\nolimits}

\newcommand{\codim}{\mathop{{\rm{codim}}}\nolimits}

\newcommand{\cl}{{\rm{closure}}}

\newcommand{\ity}{{\infty}}

\newcommand{\CL}{{\rm{CL}}}

\newcommand{\Eu}{{\rm{Eu}}}
\newcommand{\cone}{{\rm{Cone}}}
\newcommand{\NMD}{{\rm{NMD}}}
\newcommand{\fin}{\hspace*{\fill}$\Box$\vspace*{2mm}}

\newcommand{\cA}{{\cal A}}

\newcommand{\cH}{{\cal H}}

\newcommand{\cN}{{\cal N}}

\newcommand{\bC}{{\mathbb C}}

\newcommand{\bv}{{\bf{v}}}

\begin{document}

\title[Euler data for affine varieties]
 {Duality of Euler data for affine varieties}

\author{\sc Mihai Tib\u ar}

\address{Math\' ematiques, UMR 8524 CNRS,
Universit\'e de Lille 1, \  59655 Villeneuve d'Ascq, France.}

\email{tibar@math.univ-lille1.fr}







\begin{abstract} 
 We compare the Euler-Poincar\'e characteristic to the global Euler obstruction, 
 in case of singular affine varieties, and point out a certain duality among their
 expressions in terms of strata of a Whitney stratification.

\end{abstract}

\maketitle

\setcounter{section}{0}


%


  The local Euler obstruction was defined by 
MacPherson \cite{MP}, as a key ingredient for introducing Chern classes
 for singular spaces. 
Results on the local Euler obstruction have been obtained during the time by,
among others, 
A. Dubson, M.-H. Schwartz, J.-P. Brasselet, G. Gonzalez-Sprinberg, B. Teissier,
L\^e D.T, J. Sch\"urmann, J. Seade.
Some of them are surveyed in \cite{Br} and \cite{Sch-book}. For more recent 
results and generalizations one can look up 
\cite{BLS, BMPS, Sch, STV, STV-milnor}.

 For a connected singular algebraic closed affine space $Y\subset \bC^{N}$ 
we have defined in 
\cite{STV} a {\em global Euler obstruction} $\Eu(Y)$. The definition in the 
 global setting can be traced back to Dubson's viewpoint \cite{Du}. 
It immediately follows that,
for a {\em non-singular} $Y$, $\Eu(Y)$ equals the Euler 
characteristic $\chi(Y)$. The natural question that we address here
is how these two ``Euler
data'' compare to each other whenever $Y$ is {\em singular}.

Both objects, $\Eu$ and $\chi$, can be viewed as constructible functions 
with respect to some Whitney (b)-regular algebraic stratification of $Y$. 
Let us fix  such a stratification $\cA = \{
\cA_i\}_{i\in \Lambda}$ on $Y$. 
 We first show how $\Eu(Y)$ and $\chi(Y)$ can be expressed in terms of strata
such that the formulas are, in a certain sense, dual:

\begin{equation}\label{eq:1}
\Eu(Y) = \sum_{i\in \Lambda} \chi(\cA_i)\; \Eu_Y(\cA_i),
\end{equation}
\begin{equation}\label{eq:2}
\chi(Y) =  \sum_{i \in \Lambda} \Eu(\cA_i)\; \chi(\NMD(\cA_i)).
\end{equation}

  The duality consists in the observation that the formulas are obtained one
  from another
 by interchanging $\Eu$ with $\chi$. To the Euler characteristic $\chi(\cA_i)$
  of some stratum $\cA_i$ in formula (\ref{eq:1}) corresponds the global Euler obstruction
  $\Eu(\cA_i)$ of the same stratum in formula (\ref{eq:2}). The latter has the
  following meaning:  as it will be
  explained in \S \ref{glo}, the Euler obstruction  $\Eu(\bar\cA_i)$ of the
  algebraic closure $\bar\cA_i$ of $\cA_i$ in $\bC^N$ is well defined and
  depends only on the open part $\cA_i$. We may therefore set
  $\Eu(\cA_i) :=  \Eu(\bar\cA_i)$. In case of a point-stratum $\{ y \}$, we
  set $\Eu(\{ y \})= 1$.

Let us explain how the
  ``normal Euler data'' $\chi(\NMD(\cA_i))$ and
  $\Eu_Y(\cA_i)$ fit into this correspondence.
Both data are attached to a general slice $\cN_i$ of
complementary dimension of the stratum $\cA_i$ at some point $p_i\in\cA_i$.

Firstly,  $\NMD(\cA_i)$ stands for the {\em normal Morse data} of the
  stratum $\cA_i$ (after
  Goresky-MacPherson's \cite{GM}), i.e. the Morse data of $(\cN_i, p_i)$, 
see \S \ref{dual}. 

 Secondly, $\Eu_Y(\cA_i)$ denotes the {\em normal Euler obstruction}  of the
  stratum $\cA_i$, i.e. the local Euler obstruction of $\cN_i$ at $p_i$.

It is known that both data are independent on the choices of 
$\cN_i$ and of $p_i$. We refer to
  \S \ref{dual} for the definitions and more details.

We finally consider the case when $Y$ is a locally complete intersection with
  arbitrary singularities.
  We show (Proposition \ref{diff}) how the difference $\chi(Y)- \Eu(Y)$ can be
 expressed in terms of Betti
numbers of complex links and the polar invariants $\alpha_Y$ defined in \S
  \ref{glo}. If the singularities are isolated then the difference
 $\chi(Y)- \Eu(Y)$ measures the total ``quantity 
of slice-singularities'' of $Y$, see (\ref{eq:chi-iso}).

For another comparison of the Euler characteristic, namely to the
total curvature, in case of an affine hypersurface, we send the reader to \cite{ST}.

\section{Global Euler obstruction}\label{glo}

Since $Y\subset \bC^{N}$ is affine, one has a well defined {\em link at infinity} of $Y$,
denoted by $K_\ity(Y) := Y \cap S_R$. It follows from Milnor's finiteness
argument \cite[Cor. 2.8]{Mi} and from standard isotopy arguments
that $K_\ity(Y)$ does not depend on the radius $R$, provided that $R$ is large
enough.
 
Let $\widetilde Y = \cl \{ (x,T_x Y_\reg)\mid x\in Y_\reg\}\subset Y\times
 G(d,N)$ be the Nash blow-up of $Y$,
  where $G(d,N)$ is the Grassmannian of complex $d$-planes in $\bC^N $.
 Let $\nu :\widetilde Y \to Y$ denote the natural
 projection and let $\widetilde T$ denote the restriction over $\widetilde Y$ of the
 bundle $\bC^N \times U(d,N) \to \bC^N \times G(d,N)$, where $U(d,N)$
 is the tautological bundle over $G(d,N)$.  This is the ``Nash bundle'' over
 $\widetilde Y$.  We next consider a continuous,
 stratified vector field $\bv$ on a subset $V\subset Y$. The
 restriction of $\bv$ to $V$ has a well-defined canonical lifting
 $\tilde \bv$ to $\nu^{-1}(V)$ as a section of the Nash bundle
 $\widetilde T \to \widetilde Y$ (see e.g. \cite{BS}, Prop.9.1). 

We refer to \cite{STV} for other details 
concerning the following definition (which can be traced back to Dubson's
approach), and
in particular for the discussion on the independence on the choices: 
\begin{definition} \label{d:euler} Let $\tilde \bv$ be the lifting to
a section of the Nash bundle $\widetilde T$ of a stratified
vector field $\bv$ over $K_\ity(Y)= Y \cap S_R$, which is radial with respect
to the sphere $S_R$. The obstruction to extend $\tilde \bv$ as a nowhere
zero section of $\widetilde T$ within $\nu^{-1}(Y\cap B_R)$ is 
a relative cohomology class $o(\tilde \bv) \in H^{2d}(\nu^{-1}(Y \cap B_R),
\nu^{-1}(Y \cap S_R))\simeq H^{2d}(\tilde Y)$.

 One calls {\em global Euler
obstruction of} $Y$, and denotes it by $\Eu(Y)$,
 the evaluation of ${ o} (\tilde \bv)$ on the fundamental
class of the pair $(\nu^{-1}(Y\cap B_R), \nu^{-1}(Y\cap S_R))$. 

\end{definition}
 
By obstruction theory,  $\Eu(Y)$ is an integer and does not depend on the
 radius of the sphere defining the link at infinity $K_\ity(Y)$. We have shown in 
\cite[Theorem 3.4]{STV} that $\Eu(Y)$ can be expressed
in terms of polar multiplicities as follows, denoting $d = \dim Y$:
\begin{equation}\label{eq:main}
\Eu(Y) =  \sum_{j=1}^{d +1} (-1)^{d-j+1} \alpha_Y^{(j)},
\end{equation}
where:
\begin{equation}\label{eq:eulerobstr}
\alpha_Y^{(1)} := \mbox{the number of Morse points of a global generic linear
 function on } Y_\reg.
\end{equation}
After taking a general hyperplane slice $H\cap Y$, the second number is
$\alpha_Y^{(2)} := \alpha_{H\cap Y}^{(1)}$.
This continues by induction and yields a sequence of non-negative integers:
\[ \alpha_Y^{(1)}, \alpha_Y^{(2)},\cdots , \alpha_Y^{(d)},\]
which we complete by $\alpha_Y^{(d+1)}:=$ the number of points
of the intersection of $Y_\reg$ with a global generic codimension $d$ plane in
$\bC^{N}$.

 Of course $\alpha_Y^{(k)}$ depends on the embedding of $Y$ into $\bC^N$.
 Nevertheless, these invariants (and therefore, by the equality
 (\ref{eq:main}), $\Eu(Y)$ too) depend only
 on some Zariski open part of $Y$. Now, for a stratum $\cA_i$
from the stratification  $\cA = \{
\cA_i\}_{i\in \Lambda}$ of $Y$, 
the global Euler obstruction $\Eu(\bar \cA_i)$ of its Zariski closure $\bar
\cA_i$ is well-defined. However, since we have seen that this depends only on the open
part $\cA_i$, we can use the notation $\Eu(\cA_i)$ for $\Eu(\bar \cA_i)$.
This convention explains the occurrence of $\Eu(\cA_i)$ instead of
$\Eu(\bar\cA_i)$ in formula (\ref{eq:2}).

If the highest dimensional stratum is denoted by
 $\cA_0$, then we
have $\bar\cA_0 = Y$ and therefore $\Eu(Y)= \Eu(\cA_0)$.

\section{The dual formula}\label{dual}

The equality (\ref{eq:1}) was explained in \cite{STV}. 
It follows by Dubson's \cite[Theorem 1]{Du} applied to our
 setting. In case of germs of spaces a similar formula was proved in
 \cite[Theorem 3.1]{BLS} by using the Lefschetz slicing method. A different proof may be
 derived from \cite[Theorem 4.1]{BS}. For a
 more general proof, in terms of constructible functions, we send to
 \cite[(5.65)]{Sch-book}.

\smallskip

We  now give a proof of the equality (\ref{eq:2}). This can be viewed
  as a global index theorem,  similar 
to Kashiwara's local index theorem (see for this
  \cite[(5.38),(5.38)]{Sch-book}). Our proof will only use
the equality (\ref{eq:main}).

\begin{definition}\label{d:cl} (cf \cite{GM})\\
The {\em complex link} of a space germ $(X,x)$ is the general fibre in the local
Milnor-L\^e fibration defined by a general (linear) function germ at $x$.
  Up to homotopy type, this does not depend on the stratification or the
  choices of the representatives of the space or of the general function.
\end{definition}
 Let $\CL_Y(\cA_i)$ denote the {\em complex link of the stratum
 $\cA_i$} of $Y$. This is by definition the complex link of the germ $(\cN_i,
 p_i)$, where $\cN_i$ is a generic slice of $Y$ at some $p_i\in \cA_i$, of
codimension equal to the dimension of $\cA_i$. Let us remark that the complex link of a
 point-stratum $\{ y \}$ is precisely the complex link of the germ $(Y,y)$.

Let $\cone(\CL_Y(\cA_i))$ denote the cone over this complex link. We denote by
 $\NMD(\cA_i)$ the {\em normal Morse
data} at some point of  $\cA_i$,  that  is  the  pair  of  spaces 
 $(\cone(\CL_Y(\cA_i)), \CL_Y(\cA_i))$. After Goresky and MacPherson
 \cite{GM}, the local normal Morse data are local invariants up to
 homotopy and do not depend on the various choices in cause.
The complex link of the highest dimensional stratum $\cA_0$ is empty,
and we set by definition $\chi(\NMD(\cA_0)) =1$. In the same case, for the
 normal Euler obstruction we have $\Eu_Y(\cA_0) = 1$ by definition. 

\begin{theorem}\label{t:chi}
Let $Y\subset \bC^{N}$ be an algebraic closed affine space and let $\cA = \{
\cA_i\}_{i\in \Lambda}$ be some Whitney
stratification of $Y$.  Then:
\begin{equation}\label{eq:chi0}
 \chi(Y) =  \sum_{i \in \Lambda}\Eu(\cA_i)\; \chi(\NMD(\cA_i)).
\end{equation}
\end{theorem}
\begin{proof}
Take an affine Lefschetz pencil of hyperplanes in $\bC^{N}$ defined by a
linear function $l_H
:\bC^{N} \to \bC$. By the genericity of the pencil, there are only finitely
 many stratified Morse singularities of the
pencil, each one contained in a different slice. 
By the definition (\ref{eq:eulerobstr}), the
number of stratified Morse points on a stratum $\cA_i$ of dimension $>0$
 is precisely $\alpha_{\bar \cA_i}^{(\dim
\cA_i)}$.

 According to the Lefschetz
slicing method applied to singular spaces (see e.g. \cite{GM}), the space
 $Y$ is obtained
from a generic hyperplane slice $Y\cap \cH$ of the pencil, to which are
attached  cones over the complex links of each singularity of the pencil.
Goresky and MacPherson have proved that the Milnor data of a stratified 
Morse function germ is
the $(\dim \cA_i)$-times suspension of $\NMD(\cA_i)$.
At the level of Euler
characteristic, we then have:
\begin{equation}\label{eq:chi}
 \chi(Y) = \chi(Y\cap \cH) + \sum_{i \in \Lambda}(-1)^{\dim \cA_i}
 \alpha_{\bar \cA_i}^{(1)}\;  \chi(\NMD(\cA_i)),
\end{equation}
The sign $(-1)^{\dim \cA_i}$ is due to the repeated suspension of the normal
Morse data. By convention, for 0 dimensional strata $\cA_i$ we put  
$\alpha_{\bar \cA_i}^{(1)}
:= 1$, and therefore $\Eu(\bar \cA_i) =1$.
We apply formula (\ref{eq:chi}) to $Y\cap \cH$ and to the
successive generic slicings in decreasing dimensions. 
In the resulting equality, we get the sum of all
the coeffients of  $\chi(\NMD(\cA_i))$,  for each $i\in \Lambda$. We may then
 identify this sum to $\Eu(\bar \cA_i)$ via the
formula (\ref{eq:main}). This ends our proof.
\end{proof}

\section{Case of locally complete intersections}\label{compl-int}

 We consider here the case of a locally complete intersection
 $Y\subset \bC^{N}$ of dimension $d$, with arbitrary singularities. Being a
locally complete intersection implies however that the
complex link of any stratum $\cA_i$ is homotopy
 equivalent to a bouquet of spheres of dimension equal to $\codim_Y \cA_i -1$, by
 L\^e's result \cite{Le}. Let $b_{d-\dim \cA_i -1}(\CL_Y(\cA_i))$ denote the
 Betti number of this complex link.
One can then write the formula (\ref{eq:chi}) in the following
 form:
\begin{equation}\label{eq:chihyp}
 \chi(Y) = \chi(Y\cap \cH) + (-1)^d (\alpha_Y^{(1)} + \beta_Y^{(1)})
\end{equation}
where $\beta_Y^{(1)}$ collects the contributions from all the lower
 dimensional strata in the sum (\ref{eq:chi}), more precisely, under our
 assumption we have:
\[  \beta_Y^{(1)} := \sum_{i \in \Lambda\setminus \{ 0\}}
 \alpha_{\bar \cA_i}^{(1)}\;  b_{d-\dim \cA_i -1}(\CL_Y(\cA_i)).\]
According to their definitions, $\alpha_Y^{(1)}$ and  $\beta_Y^{(1)}$ are 
both non-negative integers. Their sum represents the number of $d$-cells which 
have to be attached to $Y\cap \cH$ in order to obtain $Y$.

Let us define $\beta_Y^{(k)}$ for $k\ge 2$, by:
\[ \beta_Y^{(2)} :=  \beta_{Y\cap \cH}^{(1)}\] 
and so on by induction,  for successive
slices of $Y$, as in case of the $\alpha_Y^{(k)}$-series defined
before.\footnote{We send to \cite[\S 6]{ST} for examples where the 
integers $\beta_Y^{(k)}$ are computed (but beware that we use a different
 convention for the indices $k$).} 

After repeatedly applying
(\ref{eq:chihyp}), and then using (\ref{eq:main}), we get the following 
expression of the difference among the two Euler data:
\begin{proposition}\label{diff}
\begin{equation}\label{eq:chihyp2}
 \chi(Y) - \Eu(Y)  = \sum_{k=1}^d (-1)^{d-k+1}\beta_Y^{(k)}.
\end{equation}\fin
\end{proposition} 
\begin{remark}\label{r:chi}
Let us see what becomes this difference in case $Y$ is a 
hypersurface, or a locally complete intersection, with {\em isolated singularities}.
For an isolated singular point $q\in Y$,
let $\mu_q^{\langle d-1\rangle}(Y)$ denote the Milnor
number of the local complete intersection $(Y\cap \cH , q)$ which is the result of
slicing $Y$ by a generic hyperplane $\cH$.
In case $Y$ is a hypersurface, this is the second highest Milnor-Teissier
 number in the sequence $\mu_q^{*}(Y)$.
We get:
\begin{equation}\label{eq:chi-iso}
 \chi(Y) - \Eu(Y)  = (-1)^{d}\sum_{q\in \Sing
  Y}\mu_q^{\langle d-1\rangle}(Y).
\end{equation}

Since by convention $\alpha^{(1)}_{\{
 q\}}=1$, and since $b_{d-1}(\CL_Y(\{ q\}))= \mu_q^{\langle d-1\rangle}(Y)$,
 formula (\ref{eq:chi-iso}) is indeed a particular case of 
formula (\ref{eq:chihyp2}). 
 This can be also proved by using the local Euler obstruction formula 
\cite[Theorem 3.1]{BLS}.
\end{remark}

\bigskip

{\sc Acknowledgement}  We thank J\" org Sch\" urmann for several discussions
 concerning Dubson's paper \cite{Du}, index theorems,
and his work \cite{Sch-book}.


\end{document}